\documentclass{article}
\usepackage{amssymb,amsmath, amsthm}
\usepackage{bm}
\newtheorem{theo}{Theorem}
\newtheorem{lem}[theo]{Lemma}

\title{On sums of generalized Ramanujan sums}
\author{Yusuke Fujisawa}

\begin{document}

\maketitle

\begin{abstract}
Ramanujan sums have been studied and generalized by several authors. 
For example,  Nowak \cite{nowak} studied these sums over quadratic number fields, 
and Grytczuk \cite{grytczuk} defined that on semigroups.
In this note, we deduce some properties on sums of generalized Ramanujan sums
and give examples on number fields. In particular, we have  a relational expression 
between Ramanujan sums and residues of Dedekind zeta functions.
\end{abstract}

\section{Introduction}

For  positive integers $m$ and $k$ 
the Ramanujan sum $c_k(m)$ is defined as
\begin{align*}
c_k(m)=\sum_{h \mod k \atop (h, k)=1} \exp \left( 2\pi i \frac{mh}{k} \right)
= \sum_{d \mid m, k}d \mu \left( \frac{k}{d} \right)
\end{align*}
where $\mu$ is the M{\"o}bius function. 
This sum was generalized by several authors. 
(For example, see \cite{apostol}, \cite{kiuchi}, \cite{nowak}, \cite{grytczuk}, and so on.) 
In this paper, we define generalized Ramanujan sums in another way 
and show some properties on them.

Suppose that $X$ is a non-empty set and $F_X$ is 
the set of all mappings $A:X \rightarrow \mathbb{Z}$ 
such that there are only finitely many points $x \in X$ 
such that $A(x) \neq 0$. 
We see that $F_X$ is an abelian group with respect to addition. 
For $A, B \in F_X$, we denote $A \leq B$ if $A(x) \leq B(x)$ 
for every $x \in X$. 
Let  $I_X=\{ A \in F_X : A \geq 0 \}$. 
When $X$ is the set of all prime ideals of some Dedekind domain $O$, 
we regard $I_X$ as the set of all non-zero ideals of $O$. 
 Now fix a real-valued function 
 $\mathcal{N}:I_X \rightarrow \mathbb{Z}_{> 0}$ such that 
 $\mathcal{N}(0)=1$, $\mathcal{N}(A)>1$ if $A\neq 0$, and 
$\mathcal{N}(A+B)=\mathcal{N}(A)\mathcal{N}(B)$ for all $A, B \in I_X$. 
The M{\" o}bius function $\mu$ on $I_X$ is defined as 
$\mu(A)=(-1)^{\sum_{x \in X} A(x)}$ when $A(X) \subset \{ 0, 1\}$ and 
$\mu(A)=0$ otherwise. 
For $M, K \in I_X$, we put 
\begin{align*} 
C_K(M)=\sum_{D \in I_X \atop D \leq M, K}
\mathcal{N}(D) \mu(K-D). 
\end{align*}

There are many expressions on Ramanujan sums. 
For example, 
\begin{align*}
\sum_{d|k}c_k(d)=k\prod_{p \mid k}(1-2/p), 
\end{align*}
where the product is over all prime divisors $p$ of $k$, and 
\begin{align*}
\sum_{d \mid n}c_d(m)= 
\begin{cases}
n \;\;\; \textrm{if}\; n | m, \\
0 \;\;\; \textrm{otherwise}. 
\end{cases}
\end{align*}
It is also known that 
\begin{align*}
\sum_{m=1}^{\infty}\frac{c_k(m)}{m}=-\Lambda(k) \;\;\; \textrm{if $k \neq 1$}
\end{align*}
where $\Lambda$ is the von Mangoldt function. 
Firstly, we shall show these analogues. 
Put $[x]=\sharp\{ A \in I_X : \mathcal{N}(A) \leq x \}$ for a real number $x>0$ 
when $X$ is at most countable. 
We shall show the next theorem. 

\begin{theo}\label{mt1}

\begin{itemize}
\item[(1)] 
For $K \in I_X$, we have 
\begin{align*}
\sum_{D \in I_X \atop D \leq K}C_K(D)
=\mathcal{N}(K)\prod_{p}\left( 1- \frac{2}{\mathcal{N}(A_p)} \right)
\end{align*}
where the product is over points $p \in X$ such that $K(p) \neq 0$ and 
$A_p$ is the map such that $A_p(p)=1$ and $A_p(q)=0$ if $p \neq q$.
\item[(2)] 
For $M, N \in I_X$, we have 
\begin{align*}
\sum_{D \in I_X \atop D \leq N}C_D(M)= 
\begin{cases}
\mathcal{N}(N) \;\;\; \textrm{if}\; N \leq M, \\
0 \;\;\; \textrm{otherwise}.
\end{cases}
\end{align*}
\item[(3)]
Suppose that $X$ is at most countable and  $[x]=cx+O(x^{\alpha})$ for some
 $c > 0$ and $\alpha \in [0, 1)$.  
For $K \neq 0 \in I_X$, we have 
\begin{align*}	
\sum_{M \in I_X}\frac{C_K(M)}{N(M)}=-c\Lambda(K).
\end{align*}
where $\Lambda(A)=\sum_{D \leq A}\mu(A-D)\log \mathcal{N}(D)$.
\end{itemize}
\end{theo}

Chan and Kumchev \cite{chan} studied the sums  
\begin{align*}
\sum_{m \leq x}\left(\sum_{k\leq y}c_k(m)\right)^n
\end{align*}
where $n$ is a positive integer, $x$ and $y$ are large real numbers. 
In particular, they obtain 
\begin{align*}
\sum_{m \leq x \atop k \leq y}c_k(m)=x+O(y^2). 
\end{align*}
We shall show an analogue of this expression.

\begin{theo}\label{mt2}
Suppose that $X$ is at most countable and $[x]=cx+O(x^{\alpha})$ for some
 $c > 0$ and $\alpha \in [0, 1)$. 
 
\begin{itemize} 
\item[(1)] If we fix $K \in I_X$, then
\begin{align*}
\sum_{\mathcal{N}(M) \leq x}C_K(M)
=
\begin{cases}
cx+O(x^\alpha) \;\;\; \textrm{when $K=0$,} \\
O(x^\alpha)  \;\;\; \textrm{otherwise.} 
\end{cases}
\end{align*}

\item[(2)] Put $S(x,y):=\sum_{\mathcal{N}(M)\leq x \atop \mathcal{N}(K) \leq y}C_K(M)$. 
For any $\lambda>\frac{2-\alpha}{1-\alpha}$, 
under the condition $y^{\lambda} \ll x$, considering $x \rightarrow \infty$, 
\begin{align*}
S(x, y)= cx+ o(x). 
\end{align*}
\end{itemize}
\end{theo}

\section{Preliminary}
In this section, we review or construct some basic facts of 
 arithmetical functions in a generalized situation. 
(See \cite{apostol}, \cite{cohen}, or \cite{titchmarsh}.)
Put 
$\mathcal{A}:=\{ f:I_X \rightarrow R \}$ 
where $R$ is a commutative ring. 
When $X$ is the set of prime numbers and $R \subset \mathbb{C}$, 
we may regard an elements of $\mathcal{A}$ as an arithmetical function 
in the usual case. 
Let $f$ and $g \in \mathcal{A}$. 
The Dirichelet convolution $f*g$  is 
defined as 
\begin{align*}
f*g(A)=\sum_{D \in I_X \atop D \leq A}f(D)g(A-D)
=\sum_{B, C \in I_X \atop B+C=A}f(B)g(C)
\end{align*} 
for $A \in I_X$. 
The operator $*$ on $\mathcal{A}$ is commutative, and associative. 
The identity element is the function $\delta$ such that $\delta(0)=1$ 
and $\delta(A)=0$ when $A \neq 0$. 
A function $f \in \mathcal{A}$ is invertible if and only if $f(0) \in R^{\times}$.
For simplicity, we suppose that $R=\mathbb{R}$ or $\mathbb{C}$. 
The function $\mu$ is the inverse of the function $1$ 
such that $1(A)=1$ for all $A \in I_X$, that is, $\mu*1=\delta$. 
One can see that $f=g*1$ if and only if $g=f*\mu$.  

The partial summation formula is generalized as follows. 
\begin{lem}\label{partialsum}
Suppose $[x] < \infty$ for all $x >0$.
Let $F:[1, \infty) \rightarrow \mathbb{C}$ be a $C^1$ function and $x \geq 1$. 
For $g \in \mathcal{A}$, put $S(x)=\sum_{\mathcal{N(A)}\leq x}g(A)$. 
Then 
\begin{align*}
\sum_{\mathcal{N}(A) \leq x}g(A)F(\mathcal{N}(A))
=S(x)F(x)-\int_{1}^{x}S(t)F'(t)dt.
\end{align*}
\end{lem}
This lemma is shown by the ordinary partial summation formula. So 
we omit the proof.

For a complex number $s=\sigma+it$, set  
\begin{align*}
Z(s)=\sum_{A \in I_X}\frac{1}{\mathcal{N}(A)^s}.
\end{align*}
Suppose $[x]=cx+R(x)$ and $R(x)=O(x^{\alpha})$ for some $c>0$ and $\alpha \in [0, 1)$. 
For $x>0$ and $\sigma>1$, using Lemma \ref{partialsum}, we see that  
\begin{align*}
\left| \sum_{\mathcal{N}(A)\leq x }\frac{1}{\mathcal{N}(A)^s} \right|
=\sum_{\mathcal{N}(A)\leq x}\frac{1}{\mathcal{N}(A)^{\sigma}}
=\frac{1}{x^{\sigma}}[x]+\sigma \int_{1}^{x}\frac{[t]}{t^{\sigma+1}}dt.  
\end{align*}
By our assumption, the above expression is 
\begin{align*}
\frac{c}{x^{\sigma-1}}+O(x^{\alpha-\sigma})
+\frac{\sigma cx^{1-\sigma}}{1-\sigma}
+\frac{\sigma c}{\sigma-1}+O\left(\sigma\int_{1}^{x}t^{\alpha-\sigma-1}dt \right). 
\end{align*}
Considering $x \rightarrow \infty$, we have 
\begin{align*}
 Z(\sigma) 
=\frac{\sigma c}{\sigma-1}+
 O\left( \sigma\int_{1}^{\infty}t^{\alpha-\sigma-1}dt \right).
\end{align*}
Therefore, the series $Z(s)$ is absolutely convergent for $\sigma>1$. 
Moreover,  
we can see $Z(s)$ has an analytic continuation to $\sigma > \alpha$, 
and the residue of $Z(s)$ at $s=1$ is the constant $c$. 

\section{Proof of Theorem \ref{mt1}}
For $f, g \in \mathcal{A}$ and $M, K \in I_X$, 
define the sum $S_{f, g}(M, K)$ as
\begin{align*}
S_{f, g}(M, K)=\sum_{D \leq M, K}f(D)g(K-D).
\end{align*}
Note that $S_{\mathcal{N}, \mu}(M, K)=C_K(M)$.

For $A \in I_X$, the function $\chi_A \in \mathcal{A}$ is defined as 
$\chi_A(B)=1$ when $B \leq A$ and $\chi_A(B)=0$ otherwise.
Let $v(D)=\chi_K(D)f(D)g(K-D).$
When $K$ is fixed, we see 
$S_{f, g}(M, K)=(v*1)(M)$
and 
\begin{align*}
\sum_{D \leq N} S_{f, g}(D, K)h(N-D)
&=\sum_{D \leq N}(v*1)(D)h(N-D) \\
&=(v*1*h)(N)=\sum_{D \leq N}v(D)(1*h)(N-D) \\
&=\sum_{D \leq N}\chi_K(D)f(D)g(K-D)(1*h)(N-D).
\end{align*}
Thus, we have that 
\begin{align*}
\sum_{D \leq N} S_{f, g}(D, K)h(N-D)
=\sum_{D \leq N, K}f(D)g(K-D)(1*h)(N-D)
\end{align*}
which is an analogue of Theorem 1 and 2 in \cite{apostol}.
When $f(A)=\mathcal{N}(A)$, $g(A)=\mu(A)$, and $h(A)=1$ for all $A \in I_X$, 
the above equation is 
\begin{align*}
\sum_{D \leq N} C_K(D)
=\sum_{D \leq N, K}\mathcal{N}(D)\mu(K-D)(1*1)(N-D).  
\end{align*}
Put $K=N$. We obtain 
\begin{align*}
\sum_{D \leq K} C_K(D)
&=\sum_{D \leq K}\mathcal{N}(D)\mu(K-D)(1*1)(K-D) \\
&=\sum_{D \leq K}\mathcal{N}(K-D)\mu(D)(1*1)(D) \\  
&=\mathcal{N}(K)\sum_{D \leq K}\frac{\mu(D)(1*1)(D)}{N(D)} \\
&= \mathcal{N}(K)\prod_{p}\left( 1+ \frac{\mu(A_p)(1*1)(A_p)}{\mathcal{N}(A_p)} \right) \\
&=\mathcal{N}(K)\prod_{p}\left( 1- \frac{2}{\mathcal{N}(A_p)} \right)
\end{align*}
where the product is over points $p \in X$ such that $K(p) \neq 0$ and 
$A_p$ is the map $A_p(p)=1$ and $A_p(q)=0$ if $p \neq q$. 
Hence (1) of Theorem \ref{mt1} is proved. 

Next, in order to show (2),  
fix $M \in I_X$. Then we see $S_{f, g}(M, K)=(w*g)(K)$ 
where $w(A)=\chi_M(A)f(A)$ and 
\begin{align*}
\sum_{D \leq N}S_{f, g}(M, D)h(N-D)
&= \sum_{D \leq N}(w*g)(D)h(N-D) \\
&= (w*g*h)(N) \\
&= \sum_{D \leq N}w(D)(g*h)(N-D) \\
&= \sum_{D \leq N} \chi_{M}(D)f(D)(g*h)(N-D) \\
&= \sum_{D \leq N, M}f(D)(g*h)(N-D). 
\end{align*}
Thus we obtain 
\begin{align*}
\sum_{D \leq N}S_{f, g}(M, D)h(N-D)
= \sum_{D \leq N, M}f(D)(g*h)(N-D) 
\end{align*}
which is an analogue of Theorem 3 and 4 in \cite{apostol}. 
When $f(A)=\mathcal{N}(A)$, $g(A)=\mu(A)$, and $h(A)=1$ for all $A \in I_X$, 
the above equation is 
\begin{align*}
\sum_{D \leq N}C_D(M) &= \sum_{D \leq N, M}\mathcal{N}(D) \delta(N-D) \\
&= 
\begin{cases}
\mathcal{N}(N)   \;\;\; \textrm{if $N \leq M$,} \\
0 \;\;\; \textrm{otherwise}.
\end{cases}
\end{align*}
Hence (2)  is proved. 

To show (3) we use $Z(s)$ which is defined in the previous section.
By an argument similar to that of Titchmarsh \cite{titchmarsh} p.10, we obtain that
\begin{align*}
\sum_{M \in I_X}\frac{C_K(M)}{\mathcal{N}(M)^s}
&=\sum_{M \in I_X}\frac{1}{N(M)^s}
\sum_{D \leq M, K}\mathcal{N}(D)\mu(K-D) \\
&= \sum_{D \leq K}\mu(K-D)\mathcal{N}(D)\sum_{C \in I_X}\frac{1}{\mathcal{N}(C+D)^s} \\
&= Z(s)\phi_{1-s}(K)
\end{align*}
where $\phi_{1-s}(A)=\sum_{D \leq A}\mu(A-D)\mathcal{N}(D)^{1-s}$. 
For $A \neq 0 \in I_X$, we see  
\begin{align*}
\phi_{s}(A) &=\sum_{D \leq A}\mu(A-D)\mathcal{N}(D)^{s} \\
 &= \sum_{D \leq A}\mu(A-D)\exp(s \log \mathcal{N}(D)) \\
 &= \sum_{D \leq A}\mu(A-D)
 \left( \sum_{n=0}^{\infty}\frac{ \left(s\log \mathcal{N}(D) \right)^n}{n!} \right)  \\
&= \sum_{D \leq A}\mu(A-D)
 \left( \sum_{n=1}^{\infty}\frac{ \left(s\log \mathcal{N}(D) \right)^n}{n!} \right). 
\end{align*}
Thus, we have 
\begin{align*}
\lim_{s \rightarrow 1}\frac{\phi_{1-s}(A)}{1-s}
=\lim_{s \rightarrow 0}\frac{\phi_s(A)}{s}=\Lambda(A). 
\end{align*}
The expression (3) is proved from this and 
$\lim_{s \rightarrow 1}(s-1)Z(s)=c$. 

\section{Proof of Theorem \ref{mt2}} 
Firstly, we fix $K \in I_X$. 
Then,  
\begin{align*}
\sum_{\mathcal{N}(M) \leq x}C_K(M)
&= \sum_{\mathcal{N}(M) \leq x}\sum_{D+E=K \atop D+A=M}\mathcal{N}(D)\mu(E) \\
&= \sum_{D+E=K}\mathcal{N}(D)\mu(E)\left[ \frac{x}{\mathcal{N}(D)} \right] \\
&= \sum_{D+E=K}\mathcal{N}(D)\mu(E) 
\left(c \frac{x}{\mathcal{N}(D)} + 
R\left( \frac{x}{\mathcal{N}(D)} \right) \right)  \\
\end{align*}
where $R(x)=O(x^\alpha)$. 
By the assumption, the above expression is 
\begin{align*}
cx\sum_{D+E=K}\mu(E) + O\left( x^\alpha \sum_{D+E=K}\mathcal{N}(D)^{1-\alpha} \right).
\end{align*}
Hence, $(1)$ is shown.

Next, we shall show $(2)$. 
We have 
\begin{align*}
S(x, y) &=\sum_{\mathcal{N}(M) \leq x, \atop \mathcal{N}(K) \leq y}
\sum_{D \leq M, K}\mathcal{N}(D)\mu(K-D) \\
&=  \sum_{\mathcal{N}(D+B) \leq x \atop \mathcal{N}(D+A) \leq y}
\mathcal{N}(D)\mu(A)
 \\
&= \sum_{D, A \in I_X \atop \mathcal{N}(D+A) \leq y}\mathcal{N}(D)\mu(A)
\left[ \frac{x}{\mathcal{N}(D)} \right] .
\end{align*}
By the assumption, 
\begin{align*}
S(x, y) &= \sum_{D, A \in I_X \atop \mathcal{N}(D+A) \leq y}\mathcal{N}(D)\mu(A)
\left(  c \frac{x}{\mathcal{N}(D)} + R\left( \frac{x}{\mathcal{N}(D)} \right) \right)
\\
&= cx \sum_{D, A \in I_X \atop \mathcal{N}(D+A) \leq y}\mu(A) + \sum_{D, A \in I_X \atop \mathcal{N}(D+A) \leq y}\mathcal{N}(D)\mu(A)R\left( \frac{x}{\mathcal{N}(D)} \right) 
\end{align*}
where  $R(x)=O(x^{\alpha})$. 

Since 
\begin{align*}
\sum_{D, A \in I_X \atop \mathcal{N}(D+A) \leq y}\mu(A)=
\sum_{C \in I_X \atop \mathcal{N}(C) \leq y}\sum_{A \in I_X \atop A \leq C}\mu(A)=1, 
\end{align*} 
we obtain $S(x, y)=cx+ T(x, y)$
where 
\begin{align*}
T(x, y)=\sum_{D, A \in I_X \atop \mathcal{N}(D+A) \leq y}
\mathcal{N}(D)\mu(A)R 
\left( \frac{x}{\mathcal{N}(D)} \right).
\end{align*}
Note that
\begin{align*}
T(x, y) 
\ll \sum_{D, A \in I_X \atop \mathcal{N}(D+A) \leq y}\mathcal{N}(D) 
\left( \frac{x}{\mathcal{N}(D)} \right)^{\alpha} 
= \sum_{\mathcal{N}(A) \leq y}\sum_{\mathcal{N}(D)\leq y/\mathcal{N}(A)}
x^{\alpha}\mathcal{N}(D)^{1-\alpha} \\
\ll x^{\alpha} 
\sum_{\mathcal{N}(A) \leq y}  \left( \frac{y}{\mathcal{N}(A)}\right)^{1-\alpha}
\sum_{\mathcal{N}(D)\leq y/\mathcal{N}(A)}1 
\ll x^{\alpha}y^{2-\alpha}.
\end{align*}
Hence $S(x, y)=cx+O(x^{\alpha}y^{2-\alpha})$. 
If $\frac{2-\alpha}{1-\alpha} <\lambda$ and $y^{\lambda} \ll x$, then
$T(x, y)=o(x)$. 
Therefore, Theorem \ref{mt2} is proved.

\section{Examples}

Let $F$ be a number field of degree $d$, $O_F$ the integer ring of $F$, and 
$\mathcal{I}$ is the set of all non-zero ideals of $O_F$. 
The M{\"o}bius function $\mu:\mathcal{I} \rightarrow \mathbb{C}$ for $O_F$ is defined as 
\begin{align*}
\mu(\mathfrak{a})=
\begin{cases}
(-1)^{\omega(\mathfrak{a})} \;\;\; \textrm{$\mathfrak{a}$ is square free,} \\
0 \;\;\; \textrm{otherwise}, 
\end{cases}
\end{align*}
where $\omega(\mathfrak{a})$ is the number of distinct prime factors of $\mathfrak{a}$ 
and one can define 
the Ramanujan sum as 
\begin{align*}
C_{\mathfrak{a}}(\mathfrak{b})=\sum_{\mathfrak{d} \mid \mathfrak{a}, \mathfrak{b}}
\mathcal{N}(\mathfrak{d})\mu \left( \frac{\mathfrak{a}}{\mathfrak{d}} \right)
\end{align*}
 where $\mathfrak{a}, \mathfrak{b} \in \mathcal{I}$ and 
 $\mathcal{N}(\mathfrak{a})=[O_F:\mathfrak{a}]$. 
By (1) of Theorem \ref{mt1},  
 one have 
\begin{align*}
\sum_{ \mathfrak{d} | \mathfrak{a}}C_{\mathfrak{a}}(\mathfrak{d})
=\mathcal{N}(\mathfrak{a})\prod_{\mathfrak{p} | \mathfrak{a}}
\left( 1- \frac{2}{\mathcal{N}(\mathfrak{p})} \right)
\end{align*}
for $\mathfrak{a} \in \mathcal{I}$.
By (2) of Theorem \ref{mt1}, 
 we have 
\begin{align*}
\sum_{\mathfrak{d} \mid \mathfrak{a}}C_{\mathfrak{d}}(\mathfrak{b})= 
\begin{cases}
\mathcal{N}(\mathfrak{a}) \;\;\; \textrm{if}\; \mathfrak{a} \mid \mathfrak{b}, \\
0 \;\;\; \textrm{otherwise}
\end{cases}
\end{align*}
for $\mathfrak{a}, \mathfrak{b} \in \mathcal{I}$.

The following fact is well-known.   
\begin{lem}\label{web}
(cf. Lang\cite{lang}, Chap.VI Theorem 3, or Murty and Order \cite{murty}.)
The number of ideals of $O_F$ whose norms are less than or equal to x is 
\begin{align*}
c_F x +R_F(x)
\end{align*}
where $c_F$ is the residue of the Dedekind zeta function $\zeta_F(s)$ of $F$ at $s=1$ 
and $R_F(x)=O(x^{1-\frac{1}{d}})$
\end{lem}
It is well-known that the invariant $c_F$  in the above lemma is given by 
\begin{align*}
c_F=\frac{2^{r_1}(2\pi)^{r_2}\mathcal{R}h}{W\sqrt{D}}
\end{align*}
where $r_1$ is the number of real primes, $r_2$ is the number of complex primes, 
$\mathcal{R}$ is the regulator, $h$ is the class number, $W$ is the number of roots of unity, 
and $D$ is the absolute value of the discriminant of $F$. 
The von Mangoldt function $\Lambda$ for $F$ is the function such that 
$\Lambda(\mathfrak{a})=
\log \mathcal{N}(\mathfrak{a})$ if 
$\mathfrak{a}$ is a power of a prime ideal $\mathfrak{p}$, and 
$\Lambda(\mathfrak{a})=0$ otherwise. 
Using Lemma \ref{web} and (3) of Theorem \ref{mt1}, we have 
\begin{align*}
c_F
=
-
\frac{1}{\Lambda(\mathfrak{a}) }
\sum_{\mathfrak{b}}\frac{C_{\mathfrak{a}}(\mathfrak{b})}{\mathcal{N}(\mathfrak{b})}
\end{align*}
unless $\Lambda(\mathfrak{a}) = 0$.

In addition, 
if $\lambda>d+1$ and $y^{\lambda} \ll x$, then  
\begin{align*}
\sum_{\mathcal{N}(\mathfrak{b})\leq x \atop \mathcal{N}(\mathfrak{a}) \leq y}
C_{\mathfrak{a}}(\mathfrak{b}) = c_Fx+ o(x). 
\end{align*}
by Theorem \ref{mt2}.

{\bf Acknowledgments.} The author would like to thank 
Professor Yoshio Tanigawa for telling him the equation $\sum_{m=1}^{\infty}c_k(m)/m=-\Lambda(k)$ ($k \neq 1$).

\noindent
Yusuke Fujisawa \\
Graduate School of Mathematics \\
Nagoya University \\
Chikusa-ku, Nagoya 464-8602 \\
Japan \\


\begin{thebibliography}{9}

\bibitem{apostol}
T. M. Apostol, 
{\it Arithmetical properties of generalized Ramanujan sums}, 
Pacific J. Math. {\bf 41}, no. 2, (1972), 281-293. 



\bibitem{chan}
T. H. Chan and A.V. Kumchev, 
{\it On sums of Ramanujan sums}, 
Acta arithm. {\bf 152} (2012), 1--10.

\bibitem{cohen}
H. Cohen,  
{\it Number theory. Vol. II. Analytic and modern tools},
 Graduate Texts in Mathematics, 240. Springer, New York, 2007. 


\bibitem{grytczuk}
A. Grytczuk, 
{\it On Ramanujan sums on arithmetical semigroup}, 
Tsukuba. J. Math. {\bf 16} (1992), no. 2, 315--319.


\bibitem{kiuchi}
I. Kiuchi and Y. Tanigawa, 
{\it On arithmetic functions related to the Ramanujan sum}, 
Period. Math. Hungar. {\bf 45} (2002), no. 1--2, 87--99. 


\bibitem{lang}
S. Lang, 
{\it Algebraic number theory},
2nd edition, Graduate Texts in Mathematics, 110. Springer-Verlag, New York, 1994.


\bibitem{murty}
 R. Murty and J. V. Order, 
 {\it Counting integral ideals in a number field},
Expo. Math. {\bf 25} (2007), 53--66. 

\bibitem{nowak}
W. G. Nowak, 
{\it The average size of Ramanujan sums over quadratic number fields}, 
Arch. Math. {\bf 99} (2012), 433--442.



\bibitem{titchmarsh}
E. C. Titchmarsh, 
{\it The theory of the Riemann zeta-function},
2nd edition, revised by D. R. Heath-Brown, Oxford University Press, 1986.

\end{thebibliography}
\end{document}